\documentclass[11 pt, psamsfonts]{amsart}
\usepackage{amssymb}

\def\ignore#1{\relax}

\def\g{\mathfrak g}

\def\sl{\mathfrak sl}
\def\so{\mathfrak so}
\def\sp{\mathfrak sp}
\def\p{\mathfrak p}
\def\pb{\bar p}
\def\R{{\mathbb R}}
\def\Z{{\mathbb Z}}
\def\nat{{\mathbb N}}

\def\C{{\mathbb C}}

\def\la{\lambda}
\def\La{\Lambda}

\def\N{\mathcal N}

\def\Ca{\mathcal C}

\def\one{\mathbf 1}

\def\a{{\bf a}}
\def\b{{\bf b}}
\def\W{\mathcal W}
\def\M{\mathcal M}

\def\ignore#1{\relax}
\def\om{\omega}
\def\e{\epsilon}
\def\1{{\bf 1}}

\def\D{{\mathcal D}}
\def\End{{\rm End}}
\def\Hom{{\rm Hom}}
\def\ve{\varepsilon}

\def\rhoc{{\rho\check{}}}
\def\SS{{\tilde S}}

\calclayout

\renewenvironment{proof}[1][\proofname]{\par
  \normalfont
  \topsep6\p@\@plus6\p@ \trivlist
  \item[\hskip\labelsep\bfseries
    #1\@addpunct{.}]\ignorespaces
}{%
  \qed\endtrivlist
}

\def\th@plain{%
  \let\thmhead\thmhead@plain \let\swappedhead\swappedhead@plain
  \thm@preskip.5\baselineskip\@plus.2\baselineskip
                                    \@minus.2\baselineskip
  \thm@postskip\thm@preskip
  \itshape
\renewcommand{\labelenumi}{{(\alph{enumi})\quad}}
                        \renewcommand{\labelenumii}{{(\roman{enumii})\ }}
}
\def\th@definition{%
  \let\thmhead\thmhead@plain \let\swappedhead\swappedhead@plain
  \thm@preskip.5\baselineskip\@plus.2\baselineskip
                                    \@minus.2\baselineskip
  \thm@postskip\thm@preskip
  \upshape
}
\def\th@remark{%
  \thm@headfont{\itshape}
  \let\thmhead\thmhead@plain \let\swappedhead\swappedhead@plain
  \thm@preskip.5\baselineskip\@plus.2\baselineskip
                                    \@minus.2\baselineskip
  \thm@postskip\thm@preskip
  \upshape
}

{\theoremstyle{plain}
\newtheorem{theorem}{Theorem}[section]
}

{\theoremstyle{plain}
\newtheorem{proposition}[theorem]{Proposition}
}

{\theoremstyle{plain}
\newtheorem{corollary}[theorem]{Corollary}
}

{\theoremstyle{plain}
\newtheorem{lemma}[theorem]{Lemma}
}

{\theoremstyle{plain}

}

{\theoremstyle{definition}
\newtheorem{definition}[theorem]{Definition}
}

{\theoremstyle{definition}

}

{\theoremstyle{remark}
\newtheorem{remark}[theorem]{Remark}
}

{\theoremstyle{remark}

}

\numberwithin{equation}{section}

\renewcommand{\labelenumi}{{ \theenumi.}}
\renewcommand{\labelenumii}{{(\alph{enumii})}}

\def\RepG{{\rm Rep\ }G}

\def\la{\lambda}
\def\al{\alpha}
\def\lad{\lambda^\dagger}

\def\choose #1 #2{\begin{pmatrix}#1\\#2\end{pmatrix}}

\def\qB{Br}
\def\BnN{Br_n(N)}

\def\qBb{\overline{\qB}}
\def\bNell{b^\la_\mu(N,\ell)}
\def\LNell{\Lambda(N,\ell)}

\def\x{{\bf x}}
\def\y{{\bf y}}
\def\PP{\bar P_+}
\def\ss{\tilde s}
\def\vet{\tilde\varepsilon}
\def\tLNell{\tilde\Lambda(N,\ell)}
\def\tLa{\tilde\Lambda}
\usepackage{graphicx}
\usepackage{epstopdf}

\begin{document}

\title[Fusion symmetric spaces]
{Fusion symmetric spaces and subfactors}

\author{Hans Wenzl}
\thanks{}

\address{Department of Mathematics\\ University of California\\ San Diego,
California}

\email{hwenzl@ucsd.edu}

\begin{abstract}
We construct analogs of the embedding of orthogonal and symplectic
groups into unitary groups  in the context of fusion categories.
At least some of the resulting module categories also appear in boundary conformal field theory.
We determine when these categories are unitarizable,
and explicitly calculate the  index
and principal graph of the resulting subfactors.
\end{abstract}
\maketitle

This paper is a sequel of our previous paper \cite{Wqalg}, where we
introduced a $q$-deformation of Brauer's centralizer algebra
for orthogonal and symplectic groups; this algebra had already 
appeared more or less before in \cite{Mo}, see also discussion in \cite{Wqalg}.
It is motivated by finding a deformation of  orthogonal or symplectic
subgroups of a  unitary group which is compatible  with the
standard quantum deformation of the big group.
This has been done  before on the level of coideal
subalgebras of Hopf algebras by Letzter.
However, our categorical approach also allows us to extend this
to the level of fusion tensor categories, where we find finite
analogs of symmetric spaces related to the already mentioned groups.
Moreover, we can establish $C^*$ structures, necessary for the construction
of subfactors, in this categorical setting; 
this is not so obvious to see in the setting of co-ideal algebras.

It is well-known how one can use a subgroup $H$ of a (for simplicity here)
finite group $G$ to construct a module category of the representation 
category $\RepG$ of $G$. This module category also appears in the context of
subfactors of II$_1$ von Neumann factors as follows: Let $R$ be
the hyperfinite II$_1$ factor, and let $\N=R^G\subset \M=R^H$ be
the fixed points under outer actions of $G$ and $H$. Then 
the category of $\N-\N$ bimodules is equivalent to $\RepG$, and
the module category is given via the $\M-\N$ bimodules of the 
inclusion $\N\subset\M$; its simple objects are labeled by the irreducible
representations of $H$. In particular, an important invariant called
the principal graph of the subfactor is determined by the restriction
rules for representations from $G$ to $H$.
Important examples of subfactors were constructed from fusion
categories whose Grothendieck semirings are quotients of the ones
of semisimple Lie groups. So a natural question to ask is whether
one can perform a similar construction in this context.
More precisely, can we find restriction rules for type $A$ fusion
categories which describe a subfactor as before, and
which will approach in the classical limit the usual
restriction rules from $U(N)$ to $O(N)$.

We answer this question in the positive in this paper via a fairly elementary
construction. We show that certain semisimple quotients of the $q$-Brauer
algebras have a $C^*$ structure and contain $C^*$-quotients of
Hecke algebras of type $A$. The subfactor is then obtained as
the closure of inductive limits of such algebras.
Due to its close connection to Lie groups, we can give very explicit
general formulas for its index and its first principal graph.
Observe that the Lie algebra $\sl_N$ decomposes as an $\so_N$
module into the direct sum $\so_N\oplus \p$, where $\p$ is a simple
$\so_N$-module. Then the index can be expressed explicitly in terms 
of the weights of $\p$, see Theorem \ref{explicitindex}. As before 
in the group case, it can be interpreted as the quotient of the
dimension of the given fusion category by the sum of the
squares of  $q$-dimensions of  representations
of orthogonal or symplectic subgroups whose
labels are in the alcove of a certain affine reflection group; however
in our case, there is no corresponding tensor category for the
denominator, and the $q$-dimensions differ from the ones
of the corresponding quantum groups.
Also, the restriction rules for the corresponding bimodules of
this subfactor, the first principal graph, can be derived
from the classical restriction rules via an action of the already
mentioned affine reflection
group, similarly as it was done before for tensor product rules
for fusion categories. However, in our case, the affine reflection
comes from the highest short root of the corresponding Lie algebra
in the non-simply laced case; it is also different from the one
for fusion categories in the even dimensional orthogonal case.

In the case corresponding to $O(2)$, we obtain the 
Goodman-de la Harpe-Jones subfactors for Dynkin diagrams $D_n$.
We similarly also obtain a series of subfactors for other even-dimensional
orthogonal groups. Perhaps a little surprisingly, there does not seem
to be an analogous construction which would correspond to the
group $SO(N)$ with $N$ even. Our examples for the odd-dimensional
orthogonal group and for symplectic groups seem to be closely
related to examples constructed by Feng Xu \cite{X} and Antony Wassermann
\cite{Wa} by completely different methods.  Several explicit cases are discussed in
detail at the end of the paper, as well as their connections to other
approaches, coming from boundary conformal field theory, subfactors
and tensor categories.

The first chapter  mostly contains basic material from subfactor
theory which will be needed later. In the second chapter we review
and expand material on the $q$-Brauer algebra as defined in
\cite{Wqalg}, see also \cite{Mo}. In particular, we define $C^*$-structures
for certain quotients and use that to construct subfactors.
The third chapter is mainly concerned with the finer structure of these
subfactors, such as explicit closed formulas for the index and
calculation of the first principal graph. The same techniques would
also extend to other examples, such as the ones in \cite{X}.

$Acknowledgments:$ It is a pleasure to thank Antony Wassermann,
David Jordan, Viktor Ostrik and Feng Xu for useful references,
and Fred Goodman also for technical advice.

\section{II$_1$ factors}

\subsection{Periodic commuting squares}\label{periodic}
We will construct subfactors
using the set-up of periodic commuting squares as in \cite{WHe}.
More precisely, we assume that we have increasing sequences
of finite dimensional $C^*$ algebras
$A_1\subset A_2\subset\ ... $ and 
$B_1\subset B_2\subset\ ... $ such that $A_n\subset B_n$ for all $n\in \nat$.
Let $\Lambda_n$ resp $\tLa_n$ be labeling sets for the simple components
of $B_n$ and $A_n$ respectively.
Let $G_n$ be the inclusion matrix for $A_n\subset B_n$.
If we write a minimal idempotent $p_{\mu}\in A_{n,\mu}$ as
a sum of minimal mutually commuting idempotents of $B_n$,
then the entry $g_{\la\mu}$ of $G_n$  denotes the number of those
idempotents which are in $B_{n,\la}$.
We say that our sequences of algebras are periodic with period $d$,
if there exists an $n_0\in \nat$ such that for any $n>n_0$
we have bijections $j$ between $\Lambda_n$ and $\Lambda_{n+d}$
as well as between $\tLa_n$ and $\tLa_{n+d}$ which do not change
the inclusion matrices for $A_n\subset B_n$ as well as for $A_n\subset A_{n+1}$
and $B_n\subset B_{n+1}$. This means, in particular, that $g_{j(\la)j(\mu)}
=g_{\la\mu}$ for all $\la\in\Lambda_n$, $\mu\in \tLa_n$, $n>n_0$.

The trace functional defines inner products on the algebras $A_n$ and $B_n$
by $(b_1,b_2)=tr(b_1^*b_2)$. Let $e_{A_{n+1}}$ and $e_{B_n}$ be
the orthogonal projections onto the subspaces $A_{n+1}$ and $B_n$ of $B_{n+1}$.
Then the {\it commuting square condition} says that 
$e_{A_{n+1}} e_{B_n} =e_{A_n}=e_{B_n}      e_{A_{n+1}}$ for all $n\in\nat$.
Finally, we also note that the trace $tr$ is uniquely determined on $A_n$ and
$B_n$ by its weight vectors $\a_n$ and $\b_n$ which are defined as follows:
Let $p_\mu$ be a minimal idempotent in the simple component of $A_n$ labeled by $\mu$.
Then we define $a_{n,\mu}=tr(p_\mu)$, and $\a_n=(a_{n,\mu})_\mu$, where
$\mu$ runs through a labeling set of the simple components of $A_n$. The weight
vector $\b_n$ for $B_n$ is defined similarly. The following proposition follows
from \cite{WHe}, Theorem 1.5 (where the matrix $G=(g_{\la\mu})$ defined here
would correspond to the matrix $G^t$ in \cite{WHe}):

\begin{proposition}\label{subfindex}
Under the given conditions, we get a subfactor $\N\subset\M$ whose index
$[\M :\N]$ is equal to 
$\| \a_n\|^2/\| \b_n\|^2$ for any sufficiently large $n$. Moreover,
we have $\sum g_{\la\mu} a_{n,\la}=[\M:\N]b_{n,\mu}$.
\end{proposition}

\subsection{Special periodic algebras}\label{specperiod} In general, it can be quite hard to
determine finer invariants of the subfactors, the so-called higher relative 
commutants (or centralizers) from the generating sequence of algebras.
However, under certain circumstances, this can become quite easy.
We describe such a set-up. It is a moderate abstraction of an approach
which has already been used before by a number of authors.
The reader familiar with tensor categories and module categories
should think of the algebras $A_n=\End_{\Ca}(X^{\otimes n})$ and the algebras
$B_n=\End_\D(Y\otimes X^{\otimes n})$ for $X$ an object in 
a $C^*$ tensor category  $\Ca$, and $Y$ an object in a
 module category $\D$ over $\Ca$. 
In the following, we will make the following assumptions beyond the ones
in the previous subsection:
\begin{enumerate}
\item{} The algebras $A_n$ will be monoidal $C^*$-algebras. This means we
have canonical embeddings of $C^*$ algebras $A_m\otimes A_n\to A_{n+m}$
with multiplicativity of the trace, i.e. $tr(a_1\otimes a_2)=tr(a_1)tr(a_2)$.

\item{} We have canonical embeddings $B_m\otimes A_n\to B_{n+m}$, again
with multiplicativity of the trace.

\item{} We have the commuting square condition for the sequences of algebras
$A_n\subset B_n$ and $1\otimes A_{n-1}\subset A_n$.

\item{}There exists $d\in \nat$ and a projection $p\in A_d$ such that 
$(1_m\otimes p)A_{m+d}(1\otimes p)\cong A_m$ and
$(1_m\otimes p)B_{m+d}(1\otimes p)\cong B_m$ for all $m\in\nat$.
\end{enumerate}
Examples for this set-up will be given at the end of this section and
in Section \ref{qBraueralgebras}. Moreover, any finite depth subfactor
$\N\subset \M$
(see e.g. \cite{GHJ}, \cite{EK} for definitions) produces algebras for such a set-up
as follows: Let $\M^{\otimes n}=\M\otimes_\N \M\otimes_\N\ ...\ \otimes \M$ ($n$ factors).
Obviously, $\M^{\otimes n}$ is an $\N-\N$ as well as an $\M-\N$ bimodule. One can
check that for $A_n=\End_{\N-\N}\M^{\otimes n}\subset B_n=\End_{\M-\N}\M^{\otimes n+1}$
the axioms above are satisfied; here the embedding is defined by letting the elements of
$A_n$ act on the second to $(n+1)-st$ factor of $\M^{\otimes n+1}$.
It is also possible to define these algebras in connection of relative commutants
in the Jones tower of relative commutants (see e.g. \cite{Bi} for details).
 Recall that for
 factors $\N\subset\M$ the relative commutant (or centralizer) $\N'\cap\M$
is defined to be the set $\{ b\in\M, ab=ba$ for all $a\in \N\}$.

\begin{lemma}\label{relcomm} The subfactor $\N\subset\M$ generated from
the sequences of algebras $1_m\otimes A_n\subset B_{n+m}$
has relative commutant $B_m$. The same statement also holds with
$B_{n+m}$ and $B_m$ in the last sentence replaced by $A_{n+m}$ and $A_m$.
\end{lemma}

$Proof.$ This is essentially the proof used
for Theorem 3.7 in \cite{WHe}.  Observe that by induction on $r$ and assumption 4 above,
we also have $(1_m\otimes p^{\otimes r})X_{m+rd}(1_m\otimes p^{\otimes r})
\cong X_m$ for $X=A,B$. It follows from
 Theorem 1.6 of \cite{WHe} that the dimension of the relative
commutant $\N'\cap\M$ is at most equal to the dimension of $B_m$.
The claim follows from the fact that $B_m\otimes 1_n$ commutes with
$1_m\otimes A_n$ for all $n$.

\subsection{Bimodules and principal graphs} We calculate the first principal
graph for subfactors constructed in our set-up, using fairly elementary
methods from \cite{WHe} as well as the bimodule approach. The latter
was first used in the subfactor context by Ocneanu, see e.g. \cite{EK}. For the 
connection between bimodules and principal graphs, see \cite{Bi}
and for more details compatible with our notation, see also \cite{EW}.
While most of this section has already appeared before implicitly or
explicitly, the presentation in our set-up might be useful also in other
contexts.

Pick $k$ large enough so that
$m=kd>n_0$. Hence the inclusion matrices for $A_{rd}\subset B_{rd}$
coincide for all $r\geq k$ using the bijection of simple components as described in
Section \ref{periodic}.
Let $\La_m$ and $\tilde\La_m$ be labelling sets for the simple components of
$B_m$ and $A_m$ respectively. Let  $\N$ and $\M$ be the factors
generated by the increasing sequences of algebras $A_n$ and $B_n$ respectively,
see Prop. \ref{subfindex} or Lemma \ref{relcomm}, with the $m$ there equal 0.
Both of these factors have a subfactor $\tilde\N$ generated by the subalgebras
$1_m\otimes A_n\subset A_{n+m}\subset B_{n+m}$.
We now define for each $\la\in\tilde\La_m$ an $\N-\tilde\N$ bimdule $N_\la$ as
follows: 
It is the Hilbert space completion of $ \N p_\la$ with respect to the inner
product induced by $tr$, where $p_\la$ is a minimal idempotent in $A_{m,\la}$,
the simple component of $A_m$ labeled by $\la$,
with obvious left and right actions by $\N$ and $\tilde N$.
To ease notation, we shall often refer to it as an $\N-\N$ bimodule, using
the isomorphism between $\tilde \N$ and $\N$ given by the trace preserving maps
$a\in A_n\mapsto 1_m\otimes a\in A_{n+m}$.

Similarly, we define $\M-\tilde\N$ bimodules $M_\mu$ for any $\mu\in \La_m$
which are Hilbert space completions of $\M p_\mu$, where $p_\mu$ is a minimal
idempotent in the simple component $B_{m,\mu}$ of $B_m$.
Finally, we define the inclusion numbers $b_\mu^\la$ for elements
$\la\in\tilde\La$ and  $\mu\in\La_m$ as usual (see Section \ref{periodic}).

\begin{lemma}\label{bimodule} The bimodules $N_\la$ and $M_\mu$ are
irreducible $\N-\tilde\N$ resp $\M-\tilde\N$ bimodules. We have the decomposition
$M_\mu\cong \bigoplus_\la b^\la_\mu N_\la$ as  $\N-\tilde\N$ modules.
\end{lemma}

$Proof.$ This is well-known (see e.g. \cite{EW} for more details).
 It follows from Lemma \ref{relcomm} that the endomorphism ring of
the $\M-\tilde \N$ bimodule $\M$ is given by $B_m$. Hence the 
$\M-\tilde \N$ bimodules $M_\mu$ are simple, as $p_\mu$ was chosen to
be a minimal idempotent in $B_m$. One shows similarly that also the
$N_\la$'s are simple $\N-\tilde \N$ bimodules.

Observe
that $\dim_{\N}N_\la=tr(p_\la)$ and $\dim_{\M}M_\mu=tr(p_\mu)$
(see e.g. \cite{Jo}).
Now if $p_\la$ is a minimal idempotent in $A_m$, it follows
from the definitions that  ${\rm Ind}_\N^\M N_\la := \M p_\la$
is isomorphic as an $\M-\tilde \N$ bimodule to the direct sum
$\oplus b^\la_\mu M_\mu$. By Frobenius reciprocity, see e.g. \cite{EK}, \cite{Bi}
it follows that the module $N_\la$ appears with multiplicity $b^\la_\mu$ in
$\M_\mu$, viewed as an $\N-\tilde \N$ bimodule. Hence
the $\N-\tilde \N$ bimodule $M_\mu$ has a submodule which is isomorphic
to $\bigoplus_\la b^\la_\mu N_\la$. But as
$M_\mu$ has $\N$-dimension $[\M:\N]tr(p_\mu)$, it coincides with this
submodule, by Proposition \ref{subfindex}.

\begin{theorem}\label{pringraph} Let $\N\subset\M$ be the subfactor
generated by sequences of algebras $A_n\subset B_n$ satisfying the conditions
in Section \ref{specperiod}. Then its first principal graph is given by the
inclusion graph for $A_{kd}\subset B_{kd}$ for sufficiently large $k$.
\end{theorem}

$Proof.$ It is well-known that the first principal graph is given by
the induction-restriction graph of $\M-\N$ and $\N-\N$ bimodules
appearing in the tensor products $\M^{\otimes n}$, $n\in \nat$,
where $\M^{\otimes n}=\M\otimes_\N\M\otimes\ ...\ \otimes_\N\M$
($n$ factors), see \cite{EK}, \cite{Bi}. Obviously, this graph does
not change if we replace all $X-\N$ bimodules $H$ in this setting
by $X- q\N q$ bimodules $Hq$, for $X=\M,\N$ and $q$ a nonzero
projection in $\N$. The claim can now be shown for $q=p^{\otimes k}$
where $k$ is chosen large enough so that $kd>n_0$, using 
Lemma \ref{bimodule}.

\medskip
Recall that  many examples come from module tensor categories,
where $A_n=\End_\C(X^{\otimes n})$ and
 $B_n=\End(Y\otimes (X^{\otimes n})$ for an object $X$ in
a tensor category $\C$ and an object $Y$ in the module category
$\D$ over $\C$. In this setting, the weight vectors of our
trace are given by $a_{n,\la}=\tilde d_\la/x^n$ and 
$b_{n,\mu}=d_\mu/yx^n$ for positive quantities $d_\mu, \tilde d_\la,
x$ and $y$. Then we have

\begin{corollary}\label{localind} Assuming the conditions
for the trace weights as just given, we have
 subfactors $\N\subset \M_\mu$ with index
$[\M_\mu :\N]=d_\mu^2 [\M :\N]$,
with $\N\subset \M$ as in Theorem \ref{pringraph}.
\end{corollary}

\begin{remark} There is also a second important invariant for $\N\subset\M$,
the dual principal graph. It can be analogously defined as an induction-restriction
graph between irreducible $\M-\M$ and $\M-\N$ bimodules appearing in 
the tensor powers $\M^{\otimes n}$. Its calculation is more difficult than
the first principal graph. This is quite similar to the corresponding problem
for subfactors coming from conformal inclusions and related constructions,
see e.g \cite{X1}, \cite{BEK}, \cite{EW}.
We plan to study this problem in a future publication via suitable adaptions
of techniques in those papers.
\end{remark}

\subsection{The $GHJ$-construction} 
We give a well-known and well-studied example for our current set-up,
which was first constructed in \cite{GHJ}.
 Let $G$ be a matrix with nonnegative integer
entries and norm less than 2. It is well-known that such matrices are classified
by Coxeter graphs of type $ADE$. We assume that the columns of $G$ are
indexed by the even vertices, and the rows by the odd vertices.
We define $C^*$-algebras $B_n$ by $B_0=\C^{v_e}$, and $B_1=\oplus M_{d_j}$,
where $v_e$ is the number of even vertices, and the summands of $B_1$
are labelled by the odd vertices $j$, whose dimension $d_j$ is equal to
the number of even vertices to which $j$ is connected. The embedding
$B_o\subset B_1$ is given by the inclusion matrix $G$. 
Then we define recursively $B_{n+1}$ via Jones' basic construction \cite{Jo}
for $B_{n-1}\subset B_n$. Here the trace on $B_n$  is the unique normalized
trace whose values on minimal idempotents are given by the
Perron-Frobenius vector of $G^tG$ or $GG^t$, depending on whether
$n$ is even or odd, and the vector is normalized such that $tr(1)=1$.
Then the algebra $B_{n+1}$ is generated by $B_n$, acting on itself via
left multiplication, and the orthogonal projection $e_n$ onto the subspace
$B_{n-1}$ of $B_n$, with respect to the inner product coming from the
trace.
The algebra $A_n$ is defined to be the subalgebra generated by the identity 1
and the projections $e_i$, $1\leq i<1$. 
It is well-known  that these algebras satisfy the commuting
square condition, that they are periodic with periodicity 2, and that
the Jones projections $e_i$ satisfy the conditions of the projection $p$ in
Section \ref{specperiod}. This has already been shown in \cite{GHJ}.

\section{$q$-Brauer algebras}\label{qBraueralgebras}

\subsection{Definitions}\label{qBrdef}
 Fix $N\in \Z$ and let $[N]=(q^N-q^{-N})/(q-q^{-1})$,
where $q$ is considered to be a complex number. We denote by
$H_n(q^2)$ the Hecke algebra of type $A_{n-1}$. It is given
by generators $g_1, g_2,\ ...\ g_{n-1}$ which satisfy the 
usual braid relations and the quadratic relation $g_i^2=(q^2-1)g_i+q^2$.
The $q$-Brauer algebra $\BnN$ is the complex algebra defined 
 via generators
$g_1,g_2,\ ...\ g_{n-1}$ and $e$ and relations
\begin{enumerate}
\item[(H)] The elements $g_1,g_2,\ ...\ g_{n-1}$ satisfy the relations
of the Hecke algebra $H_n(q^2)$.
\item[(E1)] $e^2=[N]e$,
\item[(E2)] $eg_i=g_ie$ for $i>2$, $eg_1=q^2e$, $eg_2e=q^{N+1}e$ and
$eg_2^{-1}e=q^{-1-N}e$.
\item[(E3)] $g_2g_3g_1^{-1}g_2^{-1}e_{(2)}=e_{(2)}=
e_{(2)}g_2g_3g_1^{-1}g_2^{-1}$,
where $e_{(2)}=e(g_2g_3g_1^{-1}g_2^{-1})e$.
\end{enumerate}

It is easy to see that this algebra coincides with the algebra defined in
\cite{Wqalg} after substituting $q$ there by $q^2$, and $e$ there by
$q^{1-N}e$ (with the $q$ of this paper); this is also compatible with
the different definition of $[N]$ in \cite{Wqalg}. We have chosen this parametrization
as it will make it easier to define a $*$ structure on it. More precisely,
if $|q|=1$,
there exists a complex conjugate antiautomorphism $b\mapsto b^*$ on $\BnN$ defined  
by 
\begin{equation}\label{staroperation}
e^*=e,\quad g_i^*=g_i^{-1},\ 1\leq i<n.
\end{equation}
It is easy to check at the relations that this operation is well-defined.

\subsection{Molev representation}\label{Molev} We  give a representation
of our algebra $\BnN$ in $\End(V^{\otimes n})$, where $V=\C^N$.
For this we use the matrices used by Molev in \cite{Mo} for the definition
of his $q$-deformation of Brauers' centralizer algebra. His defining relations
are slightly different from ours; but Molev has informed the author that
our algebra satisfies the relations of his algebra. It turns out that also
his matrices satisfy the relations of our algebras, which we will outline here.
Let $R$ be the well-known solution of the quantum Yang-Baxter equation
for type $A$. For simplicity we will use this notation for what is often denoted
as $\check{R}$. If $E_{ij}$ are the matrix units for $n\times n$ matrices,
we define the following elements in $\End(V^{\otimes 2})$:
$$R\ =\ \sum_i q E_{ii}\otimes E_{ii}\ +\ \sum_{i\neq j} E_{ij}\otimes E_{ji}\ +\ 
\sum_{i<j} (q-q^{-1}) E_{ii}\otimes E_{jj},$$
and 
$$Q= \sum_{i,j} q^{N+1-2i}\ E_{ij}\otimes E_{ij}.$$
Moreover, if $A\in \End(V^{\otimes 2})$, we define the operator $A_i\in \End(V^{\otimes n})$
by 
$$A_i\ =\ 1_{i-1}\otimes A\otimes 1_{n-1-i},$$
where $1_k$ is the identity on $V^{\otimes k}$. Then we have the following
proposition, all of whose essential parts were already proved in \cite{Mo}.
However, the relations for our algebras are slightly different, so we give
some of the adjustments of the work in \cite{Mo} to our context below.

\begin{proposition}\label{Molevrep}
The map $g_i\mapsto qR_{n-i}$ and $e\mapsto Q_{n-1}$ defines a representation $\Phi$
of $\BnN$. It specializes to the usual representation of Brauer's centralizer algebra
in $\End(V^{\otimes n})$ for $q=1$.
\end{proposition}

$Proof.$ Most of the relations are already known or are easy to check. E.g. it
is well-known that the matrices $qR_i$ satisfy the relations of the Hecke algebra
$H_n(q^2)$. Relation $(E1)$ is checked easily, and also the relations in $(E2)$
are fairly straightforward to check.  It suffices to check $(E3)$ for $n=4$.
For this observe that by \cite{Mo}, (4.16), we have 
$$Q_3R_2R_3R_1R_2Q_3=Q_1Q_3 +q^{N+1}(q-q^{-1})Q_3(R_1+q^{-1}1),$$
in our notation. Using the relation $R_i=R_i^{-1}+(q-q^{-1})1$ for
the second and third factor of the left hand side, one derives from this
$$Q_3R_2^{-1}R_3^{-1}R_1R_2Q_3=Q_1Q_3.$$
To check relation $(E3)$, observe that
$$R_1R_2(v_i\otimes v_i\otimes v_j\otimes v_j)\ =\ R_3R_2(v_i\otimes v_i\otimes v_j\otimes v_j),$$
where $(v_i)$ is the standard basis for $\C^N=V$. One derives from this 
that $R_2^{-1}R_1^{-1}R_3R_2Q_1Q_3=Q_1Q_3$.
Moreover, the same calculations above also work with $R_i$ replaced by $R_i^{-1}$
and $Q_j$ replaced by its transpose $Q_j^T$. Hence one can show as before that
$R_2R_3R_1^{-1}R_2^{-1}Q_1^TQ_3^T=Q_1^TQ_3^T$.
Transposing this, using $R_i^T=R_i$ shows the last part of the claim.

\subsection{Quotients}\label{quotients}
We can now rephrase the main
results of \cite{Wqalg} in our notation as follows:

\begin{theorem}\label{neededresults} (a) There exists a well-defined
functional $tr$ on $\BnN$ defined inductively by $tr(g_1)=q^{N+1}/[N]$,
$tr(e)=1/[N]$ and $tr(bg_n)=tr(b)tr(g_n)$ for all $b\in\BnN$.

(b) Let $\qBb_n(N)=\BnN/I_n$, where $I_n$ is the annihilator ideal of $tr$.
Then $\qBb_n(N)$ is semisimple and the inclusion $\qBb_n(N)\subset \qBb_{n+1}(N)$
 is well-defined for all $n$.
\end{theorem}

It is possible to explicitly describe the structure of the quotients $\qBb_n=\qBb_n(N)$. To do so,
we need the following definitions for the labeling sets of simple representations.
More conceptually, 
the labeling sets $\LNell$ consist of all such diagrams $\la$ for which the quantities
$d_\mu(q)\neq 0$  for any
subdiagrams $\mu\subset\la$ including $\la$ itself,
where $q^2$ is a primitive $\ell$-th root of unity and the $d_\mu$'s are defined in
Section \ref{traceweights}.

\begin{definition}\label{diagrams}  Fix  integers $N$ and $\ell$ satisfying $1<|N|<\ell$.

\noindent
(i) The set $\tLNell$ consists of all Young diagrams with $\leq N$ rows
such that the first and $N$-th row differ by at most $\ell - N$ boxes 
for $N>0$. If $N<0$, the Young diagrams have at most $|N|$ columns,
where the first and $|N|$-th column differ by at most $\ell-|N|$ boxes.

\noindent
(ii) The set $\LNell$ consists of all Young diagrams
$\la$ with $\la_i$ boxes in the $i$-th row and $\la_j'$ boxes in
the $j$-th column which satisfy

(a) $\la_1'+\la_2'\leq N$ and $\la_1\leq (\ell-N)/2$ if
$N>0$ and $\ell-N$ even,

(b)  $\la_1'+\la_2'\leq N$ and $\la_1+\la_2\leq \ell-N$ if
$N>0$ and $\ell-N$ odd,

(c) $\la_1\leq |N|/2$ and $\la_1'+\la_2'\leq \ell-|N|$ if
$N<0$ is even,

(d) $\la_1+\la_2\leq |N|$ and $\la'_1+\la'_2\leq \ell-|N|$ if
$N<0$ is odd.

\noindent
Diagrams which miss one of these inequalities only by the quantity one are called
boundary diagrams of $\LNell$; e.g. in case (a) if $\la_1'+\la_2'=N+1$.
\end{definition}

\begin{theorem}\label{qBbstructure} (\cite{Wqalg}, Section 5) Let $q^2$ be a primitive $\ell$-th root of unity,
and let  $N$ be an integer satisfying $1<|N|<\ell$.
Then the simple components of $\qBb_n=\qBb_n(N)$ 
 are labeled by the Young diagrams in $\LNell$ with
$n,n-2, n-4, ...$ boxes. If  $V_{n,\la}$ is a simple
$\qBb_n$ module for such a diagram $\la$,
 it decomposes as a $\qB_{n-1}$ module as
\begin{equation}\label{decomposition}
V_{n,\la}\ \cong \ \bigoplus_{\mu} V_{n-1,\mu},
\end{equation}
where $\mu$ runs through diagrams 
in $\LNell$ obtained by removing or,
if $|\la|<n$, also by adding a box to $\la$.
\end{theorem}
\medskip

\subsection{Path idempotents and matrix units}\label{pathsection} We will give some details
about the proof of Theorem \ref{qBbstructure} which will also be needed
for further results.
Observe that the restriction rule \ref{decomposition}
implies that a minimal idempotent $p_\mu$ in $\qBb_{n-1,\mu}$
can be written as a sum of minimal idempotents with exactly
one in $\qBb_{n,\la}$ for each diagram $\la$
 in $\LNell$ which can be
obtained by adding or subtracting a box from $\mu$.
This inductively determines a system of minimal idempotents
and matrix units of $\qBb_n(q^N,q)$ labeled by paths
resp. pairs of paths in $\LNell$ of length $n$. Such a path is
defined  to be a sequence of Young diagrams
$(\la^{(i)})_{i=0}^n$ where $\la^{(0)}$ is the  empty Young diagram,
and $\la^{(i+1)}$ is obtained from $\la^{(i)}$ by adding or removing 
a box. It follows  from the restriction rule above that
the dimension of $V_{n,\la}$ is equal to the number of paths of length
$n$ with $\la^{(n)}=\la$, and that we can label a complete system
of matrix units for the simple component $\qBb_{n,\la}$ by pairs
of such paths. We then have the following lemma:

\begin{lemma}\label{pathlemma} For each pair of paths $t_1,t_2$
in $\Lambda(N,\infty)$ with the same endpoint we can define the
matrix unit $E_{t_1,t_2}$ as a linear combination of products of
generators over algebraic functions (rational for path
idempotents) in $q$ with poles only at roots
of unity. More precisely, the formula for $E_{t_1,t_2}$
is well-defined for $q^2$ a primitive $\ell$-th root of unity
if both $t_1$ and $t_2$ are paths in $\LNell$.
\end{lemma}

$Proof.$ This was proved in  \cite{Wqalg}, Section 5.
As the result is not explicitly stated as such, we give some details here.
One observes that  the two-sided ideal
generated by the element $\bar e\in \qBb_{n+1}$ is isomorphic
to Jones' basic construction for the algebras $\qBb_{n-1}\subset\qBb_n$
(or, strictly speaking, by certain conjugated subalgebras which are
denoted by $i_1(\qBb_n)$ and $i_2(\qBb_{n-1})$, see Section
5.2 in \cite{Wqalg}). One can then define path idempotents
and matrix units inductively as it was done in \cite{RW}, 
Theorem 1.4 using the formulas for
the weights of traces, which will also be reviewed in Section
\ref{traceweights}; this  is closely related to
what is also known in subfactor
theory as the Ocneanu-Sunder path model \cite{Su}.
 The complement of this ideal is
a quotient of the Hecke algebra $\bar H_{n+1}$
for which matrix units already were more or less defined
in \cite{WHe}, p. 366.

\begin{lemma}\label{Molevcor} 
Let $p_{[1^N]}$ be the minimal idempotent in $H_N$ 
corresponding to its one-dimensional sign representation.
Then we have 
$\bar p_{[1^N]}^{\otimes 2}\qBb_{m+2N}\bar p_{[1^N]}^{\otimes 2}
\cong \qBb_{m}$ for all $m>0$.
\end{lemma}

$Proof.$  Observe that if $p\in \qBb_{2N,\emptyset}(N)$, 
the simple component labeled by the empty Young diagram
$\emptyset$,
then it follows from the restriction rule \ref{decomposition}
(see also the equivalent version below Theorem \ref{qBbstructure})
by induction on $m$ 
that $p\qBb_{m+2N}p\cong\qBb_m$ for all $m\geq 0$.
Hence it suffices to show that $ p_{[1^N]}^{\otimes 2}$
is such an idempotent.

If $q=1$ and $N>0$, $\Phi(\BnN)$ coincides with  
the commutant of the action of the orthogonal group
$O(N)$ on $V^{\otimes n}$, which is semisimple.
Moreover, the trace $tr$ is just a multiple of the pull-back
of the natural trace on $\End(V^{\otimes n})$, so
$\Phi(\BnN)\cong \qBb_n(N)$ at $q=1$. 
As $\Phi(p_{[1^N]})$ projects onto the one-dimensional
determinant representation in $V^{\otimes N}$,
the claim follows easily in that case, using Brauer duality,
 i.e. the fact that
$\Phi(\BnN)$ is equal to the commutant of $O(N)$ on $V^{\otimes n}$
for all $n$.

We will now use the fact that we can also define $\BnN$ over
the field of rational functions $\C(q)$, see \cite{Wqalg}. It follows from Lemma
\ref{pathlemma} that we can also define the path idempotents
for $\qBb_n(N)$ over that field for paths of length $n$ in
$\Lambda(N,\infty)$. As the rank of an idempotent is an
integer, the claim follows as well for $q$ a variable, and
for $q\in \C$ not a root of unity. But as 
$p_t\bar p_{[1^N]}^{\otimes 2}p_t=0$ for any path $t$
of length $2N$ in $\LNell$ which ends in $\lambda\neq\emptyset$,
we also get the  rank 0 for $\bar p_{[1^N]}^{\otimes 2}$ at
$q^2$ a primitive $\ell$-th
root of unity in $\qBb_{2N,\lambda}(N)$. This finishes the proof
for $N>0$. The proof for the symplectic case $N<0$ even goes 
the same way.

\subsection{Weights of the trace}\label{traceweights}
Using the character formulas
of orthogonal groups, one can calculate the weights of $tr$
for the algebras $\BnN$, i.e. its values at minimal idempotents
of $\BnN$. We will need the following quantities
for a given Young diagram $\la$
\begin{equation}
d(i,j)\ = \begin{cases}
\lambda_i+\lambda_j-i-j & \text{if $i\leq j$,}\\
-\lambda_i'-\lambda_j'+i+j-2 &\text{if $i> j$.}
\end{cases}
\end{equation}
Moreover, we define
 $h(i,j)$ to be the length of the hook in the Young diagram $\la$
whose corner is the box in the $i$-th row and $j$-th column. 
We can now restate \cite{Wqalg}, Theorem 4.6 in the notations 
of this paper as follows:

\begin{theorem}\label{Brweights} 
The weights of the Markov trace $tr$ for
the Hecke algebra $\bar H_n(q^2)$ are given by
$\tilde \omega_{\lambda}=\tilde d_\la/[N]^n$,
where $|\lambda|=n$, and for $\qBb_n(N)$
they are given by $\omega_{\la ,n} = d_\la/[N]^n$, where 
\begin{align}\label{Heckeweights1}
\tilde d_\la\ =\ \prod_{(i,j)\in\la} \frac{[ N+j-i]}{[h(i,j)]},
\quad
d_\la\ =\ 
\prod_{(i,j)\in\la} \frac{[ N+d(i,j)]}{[h(i,j)]}, \notag
\end{align}
where $\la$ runs through all the Young diagrams in
$\tLNell$ with $n$ boxes for $\bar H_n(q^2)$, and through
all Young diagrams in
$\LNell$  with $n, n-2, n-4, ... $
boxes. for $\qBb_n$.
\end{theorem}

\begin{lemma}\label{weightpositivity}
The weights $\om_{\la,n}$ are positive for all $\la\in\LNell$
if and only if $q^2=e^{\pm 2\pi i/\ell}$ with $\ell>N$ and

(a) $N>0$ and $\ell-N$ even or

(b) $N<0$ odd.
\end{lemma}

$Proof.$  The weights can be rewritten for our choice of $q$ as
$$\omega_{\la ,n}\  =\ \frac{\sin^n(\pi/\ell)}{\sin^n(N\pi /\ell)^n}\
\prod_{(i,j)\in\la} \frac{ \sin (N+d(i,j))\pi/\ell}{\sin(h(i,j)\pi/\ell)}.$$
As $h(i,j)\leq h(1,1)=\la_1+\la_1'-1<\ell$ for all boxes $(i,j)$ of $\la$,
it follows that all factors in the formula above are positive for $N>0$
(negative for $N<0$) except possibly the ones in the numerator under
the product.  If $N>0$ and $\ell-N$ odd, one checks that 
for the diagram $[\ell-N+1)/2]$ we have $\om_{\la,|\la |}<0$.
By the same argument, one shows that $\om_{\la,|\la |}<0$
for $\la = [(|N|+1)/2]$ and $N<0$.
In the other two cases, one checks that $0<|d(i,j)|<\ell$ for
all boxes $(i,j)$ of a diagram $\la\in\LNell$.

\medskip

\subsection{$C^*$-quotients}
\begin{proposition}\label{Cstar} If the weights $\om_{\la, n}$
are positive for all $\la\in\LNell$,
the star operation defined by $e^*=e$ and by $g_i^*=g_i^{-1}$
makes the quotients $\qBb_n$ into $C^*$ algebras. 
\end{proposition}

$Proof.$ The proof goes by induction on $n$, with the claims for $n=1$
and $n=2$ easy to check.
 By \cite{Wqalg}, the two-sided ideal $I_{n+1}$ generated by $e$ in $\qBb_{n+1}$
is isomorphic to Jones' basic construction for $\qBb_{n-1}\subset\qBb_n$,
see also the remarks before Lemma \ref{Molevcor}.
In particular, this ideal is spanned by elements of the form
$b_1eb_2$, with $b_1,b_2\in i_1(\qBb_n)$, where $i_1(a)=\Delta_{n+1}a\Delta_{n+1}^{-1}$,
with $\Delta=(g_1g_2\ ... g_{n-1})(g_1\ ...\ g_{n-2})\ ...\ g_1$.
By induction assumption and properties of Jones' basic construction,
this ideal  has a $C^*$ structure given by
$(b_1eb_2)^*=b_2^*eb_1^*$. This coincides with the $*$ operation defined
before algebraically. It was shown in \cite{Wqalg} that $\qBb_{n+1}\cong
I_{n+1}\oplus \bar H_{n+1}$, where $\bar H_{n+1}$ is a semisimple quotient
of the Hecke algebra $H_{n+1}$ whose simple components are labeled by
the Young diagrams $\la\in\LNell$ with $n+1$ boxes. 
All these simple representations satisfy the $(k,\ell)$ condition in \cite{WHe}.
It follows from that paper that the map $g_i^*=g_i^{-1}$ induces a $C^*$
structure for any such representation. This finishes the proof.

\begin{theorem}\label{subfcon}
For each choice of $N$ and $\ell$ with $q^2=e^{\pm 2\pi i/\ell}$, and
for each nonnegative  integer $m$ we obtain 
a subfactor $\N\subset \M$ with $\N'\cap \M=\qBb_m$ and with
index
$$[\M:\N]\ =\ [N]^m\ \frac{\sum_{\mu\in\tLNell} \tilde d_\mu^2}{\sum_{\la\in\LNell}d_\la^2},$$
with notations as in Section \ref{traceweights}. Moreover,  its first principal
graph is given by the inclusion graph for $\bar H_{2Nk}\subset\qBb_{2nk+m}$
for any sufficiently large $k$.
\end{theorem}

$Proof.$ Let us first check conditions 1-4 in Section \ref{specperiod}
with $A_n=\bar H_n$ and $B_n=\qBb_n(N)$ for $q=e^{\pi i/\ell}$
and $1<|N|<\ell$. Condition 1 is well-known and was checked 
in e.g.\cite{WHe}. Similarly, Cond. 2 follows from the results in \cite{Wqalg},
using the map $b\otimes g_i\in \qBb_m\otimes\bar H_n\mapsto bg_{m+i}$.
Condition 3  means that the conditional expectation from $\qBb_{n+1}$ to
$\qBb_n$ maps $\bar H_{n+1}$ onto $\bar H_n$. But as any element of
$\bar H_{n+1}$ can be written as a linear combination of elements of the form
$ag_nb$, with $a,b\in\bar H_n$, we have for any $c\in\qBb_n$ that
$$tr(ag_nbc)=tr(g_n)tr(abc)=tr(E_{\bar H_n}(ag_nb)c).$$
Hence the commuting square condition is satisfied for any four algebras
of the type above.
Finally,  Condition 4 follows for $d=2N$ and the
projection $p=p_{[1^N]}^{\otimes 2r}$ from  Lemma \ref{Molevcor}.

The periodicity condition for $\bar H_n$ was 
 shown in \cite{WHe} by proving that $\bar p_{[1^N]}\bar H_{m+N}\bar p_{[1^N]}
\cong \bar H_{m}$. This induces an injective map $\tLNell_m\to\tLNell_{m+N}$
by adding a column of $N$ boxes to the given Young diagram which has to
become surjective for sufficiently large $m$ by definition of $\tLNell$.
The $2N$ periodicity for the algebras $\qBb_n(N)$ follows similarly
using Lemma \ref{Molevcor}, or see \cite{Wqalg}.

\section{$S$-matrix} 

We will need certain well-known identities, 
which can be found in \cite{Kc}, except for one case, which is a 
variation of the other ones. Because of this, we review the material
in more detail. This might also be useful to the non-expert reader,
as the identities needed here can be derived by completely elementary
methods.

\subsection{Lattices}
Let $M\subset L\subset \R^k$ be two lattices of full rank. This means that
they are isomorphic to $\Z^k$ as abelian groups, and each of them
spans $\R^k$  over $\R$. Moreover, we assume that we have an inner
product on $\R^k$ such that $(\x,\y)\in \Z$ for all $\x,\y\in M$.
We define the dual lattice $M^*$ to be the set of all $\y\in\R^k$ such that
$(\x,\y)\in \Z$ for all $\x\in M$; the dual lattice $L^*$ is defined similarly.
Obviously $M\subset L$ implies $L^*\subset M^* $.
Finally, we also assume that $A=L/M$ is a finite abelian group.
Then each $\gamma\in M^*$ defines a character of $A$ via
the map $e^\gamma: \x\in L\mapsto e^{2\pi i(\gamma,\x)}$.
In particular, one can identify the group dual of $A$ with $M^*/L^*$.
Define the matrix 
$\SS=\frac{1}{|L:M|^{1/2}}(e^\gamma(\x))$, 
where $\gamma$ and $\x$ 
are representatives for the cosets $M^*/L^*$ and $L/M$.
Then $\SS$ is  the character matrix of $A$ up to a multiple 
and one easily concludes that it is unitary. More precisely,
we can view it as a unitary operator between Hilbert spaces
$V$ and $V^*$ with orthonormal bases labeled by the elements 
of $L/M$ and $M^*/L^*$ respectively.

\subsection{Weights of traces} We will primarily be interested in 
lattices related to root, coroot and weight lattices of orthogonal and
symplectic groups.  We define the lattices
\begin{equation}\label{PQ}
Q=\{ \x\in \Z^k,\ 2|\sum x_i\}\quad {\rm and} \quad P=\Z^k\cup (\ve +\Z^k),
\end{equation}
where $\ve$ is the element in $\R^k$ with all
its coordinates equal to $1/2$.
Observe that $P^*=Q$ with respect to the usual scalar product of $\R^k$.
Moreover, one can identify coroot and weight lattices of $\so_{2k}$ or $\so_{2k+1}$
with $Q$ and $P$ respectively. In particular, we define for any $\gamma\in P$
the functional $e^\gamma:\R^k\to \C$ by $e^\gamma(\x)=e^{2\pi i(\gamma,\x)}$.
 The Weyl group of type $B_k$ 
acts as usual via permutations and sign changes on the coordinates.
 Let $a_W=\sum_w\ve(w)w$, where $\ve(w)$ is the sign of the element $w$.
Then the characters $\chi_\la$  for $\so_{2k+1}$ resp for $\sp_{2k}$
are given by $\chi_\la=a_w(e^{\la+\rho})/a_w(e^\rho)$,
where $\rho=(k+1/2-i)$ for $\so_{2k+1}$ and $\rho=(k+1-i)$
for $\sp_{2k}$, and $W$ is the Weyl group of type $B_k$.

We will also need the somewhat less familiar character formulas for
the full orthogonal group $O(N)$:
Recall that the irreducible representations of $O(N)$ are labeled by
Young diagrams $\la$ with at most $N$ boxes in the first two columns. 
$O(N)$-modules labeled by Young diagrams
$\la\neq\lad$ restrict to isomorphic $SO(N)$-modules if and only if $\la_1'=N-(\lad)_1'$
and $\la_i'=(\lad)_i'$ for $i>1$. Hence if $g=exp(\x)$ is an element in $SO(N)$,
 it suffices to consider the quantities
$\chi_\la(g)=\chi_\la(\x)$ for $\la$ with at most $k$ rows for $N=2k$ or $N=2k+1$.
We can now express the weights of Theorem \ref{Brweights}
in terms of these characters; in fact the formulas in Theorem \ref{Brweights} were
derived from these characters, see \cite{Ko} and \cite{Wqalg}.

\begin{lemma}\label{weightcharacters}
Let $d_\la,\ \tilde d_\la$ be as in Theorem \ref{Brweights} for $q=e^{\pi i/\ell}$.
Moreover, we define for $|N|=2k$ or $N=2k+1$ the vector
$\rhoc\in \R^k$ by $\rhoc=((|N|+1)/2-i)_i$.  By the discussion above,
it suffices to evaluate $\chi^{O(N)}(\rhoc/\ell)$ for Young diagrams $\la$
with $\la_1'\leq N/2$, which will be assumed in the following.

(a) If $N=2k+1>0$, then $d_\la=\chi_\la^{O(N)}(\rhoc/\ell) =
\chi_\la^{SO(N)}(\rhoc/\ell)$.

(b) If $N=2k>0$ and $\la_1'\leq k$, then $d_\la=m(\la)\det(\cos(l_j\rhoc_i)/\det(\cos(k-j)\rhoc_i)$,
where $l_j=(\la+\rho)_j=\la_j+k-j$ and
where $m(\la)=2$ or $1$, depending on whether $\la$ has exactly
$k$ rows or not.

(c) If $N=-2k$, then $d_\la=(-1)^{|\la |}\chi_{\la^t}(\rhoc/\ell)$ for the symplectic
character labelled by the transposed diagram $\la^T$.

(d) We have $\tilde d_\la=\chi_\la^{SU(N)}(\rho/\ell)$ for $N>0$
and  $\tilde d_\la=(-1)^{|\la|}\chi_{\la^T}^{SU(N)}(\rho/\ell)$  for $N<0$,
where
$\rho=((|N|+1)/2-i)\in \R^{|N|}$.
\end{lemma}

$Proof.$ Observe that $\rhoc$ is the element $\rho$ of the Cartan subalgebra
of $\sl_N$, viewed as an element of the Cartan subalgebra of 
the Lie subalgebra $\so_N$ or $\sp_N$, depending on the case.
The proof now goes as e.g the proof of Theorem 4.6 in \cite{Wqalg}, 
which is essentially the one of \cite{Ko}.
The fact that these arguments also work for the special quotients
$\qBb_n$ follows from the proof of \cite{Wqalg}, Theorem 5.5.

\begin{remark} Let $\Delta_+$ be the set of positive roots of a semisimple
Lie algebra and $|\Delta_+|$ be its cardinality.
As usual, we can express the Weyl denominator
in $\chi_\la(\rhoc/\ell)$ in product form as
\begin{equation}\label{denominator}
\Delta(\rhoc/\ell)\ =\ \prod_{\al >0}(e^{(\al,\rhoc)\pi i/\ell}-e^{-(\al,\rhoc)\pi i/\ell})
\ =\ (-i)^{|\Delta_+|} \ \prod_{\al>0} 2\sin ((\al,\rhoc)\pi/\ell).
\end{equation}
\end{remark}

\subsection{Usual $S$-matrices}
As usual, we pick as dominant chamber $C_+$ the regions given by
$x_1>x_2>\ ...\ x_k>0$ for Lie types $B_k$ and $C_k$.
We also choose the fundamental domains $D$ with respect to
the translation actions of $M,M^*,L,L^*$ such that it has 0 in its
center; here the lattices $M$ and $L$ will be certain multiples
of the lattices $P$, $Q$ or $\Z^k$ to be specified later.
Let $\PP$ be the intersection of $M^*$ with the fundamental
alcove  $D\cap C_+$. 

Observe that we also obtain a representation of the Weyl group
$W$ on the vector spaces $V$ and $V^*$. 
Then it is easy to check that $a_W(V^*)$ has an orthonormal  basis 
$|W|^{-1/2}a_w(e^\gamma)$, with $\gamma\in\PP$,
and we can define a similar basis $a_W(\x)$ for $a_W(V)$.
Let $S$ be the matrix which describes the action of $\SS_{|a_W(V)}$
with respect to that basis. Then it is not hard to check (and we will do 
a slightly more complicated case below) that its coefficents are given by
\begin{equation}\label{smatrixentry}
s_{\gamma,\x}=\frac{1}{|L:M|^{1/2}}\sum_w\ve(w) e^{2\pi i(w.\gamma, \x)}.
\end{equation}
If $L$ is the weight lattice of a simple Lie algebra, the entry $s_{\gamma,\x}$
is the numerator of Weyl's character formula for the dominant weight $\la=\gamma-\rho$,
up to the factor $|L:M|^{-1/2}$. As the columns of the unitary matrix $S$ have norm one,
it follows that 
\begin{equation}\label{squaresums}
\sum_\la \chi_\la^2(\x) = \frac{|L:M]}{\Delta^2(\x)},
\end{equation}
where $\Delta$ is the Weyl denominator, and the summation goes over the dominant weights
$\la$ such that $\la+\rho\in \PP$. We are now in the position to  prove
some cases of the following proposition:

\begin{proposition}\label{squaresum}
Let $\LNell_{ev}$ be the subset of $\LNell$ consisting of Young diagrams with an even
number of  boxes.
Then we have
$$\sum_{\la\in \LNell_{ev}} d_\la^2\ =\ \frac{\ell^k}{b(N)}\ \prod_{\al>0}\frac{1}{4\sin^2 (\al,\rhoc)\pi/\ell},$$
where $\rhoc=((|N|+1)/2 -i)$ and $\al>0$ runs through the positive roots of $\so_N$ for
$N>0$ and of $\sp_{|N|}$ for $N<0$ even, and $b(N)=2$ for $N=2k>0$,
and $b(N)=1$ otherwise.
\end{proposition}

$Proof.$
 Let us consider the case $N=2k+1>0$,
with $P$ and $Q$ as in \ref{PQ}.  Let $L=\ell^{-1}\Z^k$
and let $M_1=Q$ and $M_2=\Z^k$.  Then we have 
$M_1^*=P$, and $M_2^*=\Z^k$. Now  observe that $M_1^*$ is
the weight lattice of $\so_N$, and the elements $\gamma\in \PP$ are
in 1-1 correspondence with the dominant weights $\la$ of $\so_N$
satisfying $\la_1\leq (\ell-N)/2$, via the correspondence $\gamma=\la+\rho$.
Moreover, $|L:M_1|=2\ell^k$. Hence it follows from Eq \ref{squaresums}
that $\sum \chi_\la^2(\rhoc)=2\ell^k/\Delta^2(\rhoc)$.
Playing the same game for the lattice $M_2$, we now only get the sum over
the characters $\chi_\la^2$ for which $\la+\rho$ is in $\Z^k$, which is only
half as large as before. Hence also the sum over the characters $\chi_\la^2$
for which $\la\in\Z^k$ has to have the same value. 
This sum coincides with the right hand side of the statement for $N>0$ odd,
by the restriction rules for $O(N)$ to $SO(N)$ (see Lemma \ref{weightcharacters}
and its preceding discussion).

The symplectic case $N=-2k<0$ goes similarly. Here we define $M\subset L=\ell^{-1}P$,
and with $L^*=\ell Q\subset M^*=\Z^k$. Then it follows that
$\sum d_\la^2=2\ell^k/\Delta^2(\rhoc/\ell)$, where the summation goes over
all diagrams $\la$ such that $\la^T\in \LNell$. Playing the same game for $M=P$
and $M^*=Q$, we get $\sum d_\la^2=\ell^k/\Delta(\rhoc/\ell)$, where now
the summation goes over all even, or over all odd diagrams in $\LNell$, depending on
whether the sum of coordinates of $\rho=(k+1-i)$ is odd or even.
In each case, we obtain that
$\sum_{ev} d_\la^2=\ell^k/\Delta(\rhoc/\ell)$. We have proved the proposition
except for the case $N=2k>0$, for which we need a little more preparation.

\subsection{Another $S$-matrix}\label{anothers}  We now consider a slight generalization of the above.
Observe that we can define a second sign function $\vet$  for $W=W(B_k)$ which 
coincides with the usual sign function on its normal subgroup $W(D_k)$,
while we have $\vet(w)=-\ve(w)$ for $w\not\in W(D_k)$. It is easy to see that
also in this case we have $\vet(vw)=\vet(v)\vet(w)$ for all $v,w\in W$.  We define 
$\tilde a_W=\sum \vet(w) w$, and also denote the corresponding operators
on the various (quotient) lattices and on the vector spaces $V$ and $V^*$
by the same symbol. One observes that now we get an orthonormal  basis for
$\tilde a_W(V^*)$ of the form $\b_\gamma=|Stab(\gamma)|^{-1/2}|W|^{-1/2}\tilde a_W(e^\gamma)$, 
labeled by the elements of $\PP$ which now consist of the $\gamma\in D$ such that
$\gamma_1>\gamma_2>\ ...\ \gamma_k\geq 0$.
Observe that $|Stab(\gamma)|$ is equal to 1 or 2, depending on whether $\gamma_k>0$ or
$\gamma_k=0$.  One similarly defines a basis for $\tilde a_W(V)$.
 Let $\x$ be such that $Stab(\x)=1$, i.e. $x_k>0$, and let
$b_\x=|W|^{-1/2}\tilde a_W(\x)$.  Then,
writing $M^*/\ell L^*$ as a collection of $W$ orbits, we obtain

\begin{align}\notag
\SS \b_\x\ &\ =|W|^{-1/2} \sum_{\la \in\PP} \sum_{v,w\in W} \frac{1}{|Stab_W(\gamma)|}
\vet(w)\ss_{v.\gamma,w.\x} v.\gamma \notag\\
&\ =\ \sum_{\la\in\PP}\sum_{v}(\sum_w \vet(w)\ss_{w.\gamma,\x}\frac{1}{|Stab_W(\la)|})\vet(v)v.\gamma.\,
\notag
\end{align}
where we replaced $\vet(w)$ by $\vet(v)\vet(w^{-1}v)$, $\ss_{v.\la, w.\x}$ by $\ss_{w^{-1}v.\gamma,\x}$
and finally also substituted $w^{-1}v$ by $w$.
We see from this that the coefficient of $v.\gamma$ is equal to 0 if $\gamma$ has a nontrivial stabilizer
except in the case when $\gamma_k=0$.
Hence it follows that $\SS$ maps $a_W(V)$ into $a_W(V^*)$. Taking bases $(\tilde a_W(\gamma))_{\gamma\in P_+}$
and $(\tilde a_W(\x))$, we see that $\SS_{| a_W(V)}$ can be described by the matrix $S=(s_{\gamma,\x})$ 
whose coefficients are given for $\x$ with trivial stabilizer by
\begin{equation}\label{asmatrixentry}
s_{\gamma,\x}=|Stab(\gamma)|^{-1/2}|L:M|^{-1/2}\sum_w\ve(w) e^{2\pi i(w.\gamma, \x)}.
\end{equation}

\subsection{Squares of characters}\label{squares}
Using the discussion before and the formulas of Lemma \ref{weightcharacters}
 it is not hard to see that for $N$ even and $\la_1'\leq N/2$ we can write

$$\chi_\la^{O(N)}= m(\la)\tilde a_W(e^{\la+\rho})/\tilde a_W(e^\rho),$$
where $m(\la)=2$ or $1$ depending on whether $\la$ has exactly $k$ rows or not.
 In particular, applying this to the trivial representation,
we obtain $2\Delta(\rho) = \tilde a_W(e^\rho)$.

Let $P$ and $Q$ be as in \ref{PQ}, and set $L=\ell^{-1}P$ and
$M=\Z^k$. Then $L^*=\ell Q\subset M^*=\Z^k$, and it is 
 easy to see that all of these lattices  are $W=W(B_k)$-invariant. 
Moreover, let $\rhoc/\ell=(k+1/2-i)/\ell\in \ell^{-1}P=M^*$. Then it follows for $N=2k$
and $\ell$ even that 
$$\sum_{\la\in \LNell} \chi^2_\la(\rhoc\ell)\ =\ \frac{1}{\Delta^2(\rhoc\ell)}\sum_{\la_{k+1}
= 0,\la_1\leq (\ell-N)/2} (\tilde a_W(e^{\la+\rho})(\rhoc)\ell)^2
\ =\ \frac{|L:M]}{2\Delta^2(\rhoc)}\sum_\la s_{\la,\rhoc/\ell}^2.$$
Now observe that the matrix $S$ is unitary and that $[L:M]=2\ell^k$. Moreover, by e.g. 
Prop. \ref{subfindex} and Theorem \ref{subfcon}
the square sum over odd diagrams must be
equal to the square sum over even diagrams. Hence we obtain for $N>0$ even, and $\ell$ even that
\begin{equation}\label{O(2k)}
\sum_{\lambda\in\LNell_{ev}} d_\la^2\ =\ 
\frac{\ell^k}{2\Delta^2(\rhoc)},
\end{equation}
where $\LNell_{ev}$ denotes the set of diagrams in $\LNell$  with an even number of boxes.
This finishes the last case of the proof of Proposition \ref{squaresum}

\subsection{Calculation of index} As usual, identify the Cartan algebra of $\sl_N$ with
the diagonal $N\times N$ matrices with zero trace. The embedding of the Cartan algebras
of an orthogonal or symplectic subalgebra is given via diagonal matrices for
which the $N+1-i$-th entry is the negative of the $i$-th entry, for $1\leq i\leq N/2$.
Hence, if $\e_i$ is the $\sl_N$ weight given by the projection onto the $i$-th diagonal
entry, we have $(\e_{N+1-i})_{|\so_N}=(-\e_i)_{|\so_N}$, with a similar identity
also holding for symplectic subalgebras. Using our description of coroot
and weight lattices of orthogonal and symplectic Lie algebras as sublattices of $\R^k$,
and defining $\phi_i$ to be the projection onto the $i$-coordinate,
we see that  $(\e_{N+1-i})_{|\so_N}=-\phi_i=(-\e_i)_{|\so_N}$.
This allows us to describe the decomposition of $\sl_N$ as an $\so_N$ resp. $\sp_N$
module as follows: We have
 \begin{equation}\label{slNdecomposition}
\sl_N=\so_N\oplus\p,\quad {\rm resp.}\quad \sl_N=\sp_N\oplus \p
\end{equation}
where $\p$ is the nontrivial irreducible submodule in
the symmetrization of the vector representation of $\so_N$,
resp.  $\p$ is the nontrivial irreducible submodule in
the antisymmetrization of the vector representation of $\sp_N$.
The   nonzero weights $\omega>0$ of $\p$ coming from positive roots of $\sl_N$
and the multiplicity $n(\p)$  of the weight 0 in $\p$ 
are given by 

(a)\ $2\phi_i, \phi_i$ and $\phi_i\pm\phi_j$ for $1\leq i<j\leq k$
with $n(\p)=k$
for $\so_N$ with  $N=2k+1$ odd,

(b)\  $2\phi_i$ and $\phi_i\pm\phi_j$ for $1\leq i<j\leq k$ with $n(\p)=k-1$ 
for  $\so_N$ with $N=2k$ even,

(c) \  $\phi_i\pm\phi_j$ for $1\leq i<j\leq k$ with $n(\p)=k-1$ 
for  $\sp_{|N|}$ with $N=-2k<0$ even.

\begin{theorem}\label{explicitindex}
The index of the subfactor $\N\subset\M$ obtained from the inclusions of algebras
$\bar H_n(q)\subset \qBb_n(q^N,q)$ is given by
$$[\M:\N]= b(N)\ell^{n(\p)}\prod_{\omega>0}\frac{1}{4\sin^2 (\om,\rhoc)\pi/\ell},$$
where the product goes over the weights $\omega>0$ of $\p$ coming from positive
roots of $\sl_N$, as listed above, $n(\p)$ is the multiplicity of the zero weight in $\p$,
and $b(N)$ and $\rhoc$ are as in Prop. \ref{squaresum}.
\end{theorem}

\begin{corollary}\label{asymptotics} If $q=e^{\pi i/\ell}\to 1$, the index $[\M:\N]$ 
goes to $\infty$ with asymptotics $\ell^{\dim \p}$.
\end{corollary}

$Proof.$  We use Theorem \ref{subfcon}, where the denominator has been calculated 
in Proposition \ref{squaresum}. The numerator follows from a standard argument 
for $S$-matrices for Lie type $A$, see \cite{Kc}, versions of which have
also been used in this section. For an elementary calculation,
see \cite{E}.

\begin{remark}\label{indexremark} It is straightforward to adapt our index 
formula to subfactors related to other  fixed points $H=G^\alpha$  of an order two
automorphism $\alpha$  of a compact Lie group $G$, up to some
integer (or perhaps rational) constant $b(H,G)$.  Again, $\p$ would be
the $-1$ eigenspace of the induced action of $\alpha$ on the Lie algebra $\g$,
and the same $S$-matrix techniques applied in this section would go through.
E.g. our formulas for $N=3$ and $\ell$ odd coincide with the
ones at the end of \cite{X} for even level of $SU(3)$, up to a factor 3. 
This is to be expected 
as in our case only those diagrams appear in the principal graph (see
next section)  which also label
representations of the projective group $PSU(3)$.
\end{remark}

\subsection{Restriction rules and principal graph} It follows from Theorem \ref{pringraph}
that the principal graph of $\N\subset\M$ is given by the inclusion matrix for
$\bar H_{2k}\subset\qBb_{2k}$ for $k$ sufficiently large. This still leaves the question
how to explicitly calculate these graphs. Observe that in the classical case $q=1$
these would be given by the restriction rules from the unitary group
$U(N)$ to $O(N)$, for $N>0$.
Formulas for these restriction coefficients have been well-known, see
e.g. \cite{Wy} (see Theorems 7.8F and 7.9C), Littlewood's formula
(see e.g.  \cite{KT}, Section 1.5, and the whole paper for additional results).
Another approach closely related to the setting of fusion categories
can also be found in \cite{Wquot}.

Let $b^\la_\mu(N)$ be the multiplicity of the simple $O(N)$-module $V_\mu$
in the $U(N)$ module $F_\la$, for $N>0$, where $\la$, $\mu$ are Young diagrams.
It is well-known that for fixed Young diagrams $\la$ and $\mu$,
the number $b_\mu^\la(N)$ will become a constant $b_\mu^\la$
for $N$ large enough.
 Fix now also $\ell>|N|$.
 We define similar coefficients in our setting as follows:
Recall that the simple components of $\bar H_n$ are
labeled by the diagrams in $\tLNell_n$ and the ones of $\qBb_n$ by
the diagrams in $\LNell$. We then define for $\la\in\tLNell$ and
$\mu\in\LNell$ the number
$\bNell$ to be the multiplicity of a simple $\bar H_{n,\la}$ module 
in a simple $\qBb_{n,\mu}$ module.

In the following lemma the symbol $\chi_\mu$ will also be used for the $O(N)$
character corresponding to the simple representation labeled by the Young diagram
$\mu$. Moreover, we also denote by $\qBb_\infty$ the inductive limit of
the finite dimensional algebras $\qBb_n$ under their standard inclusions,
for fixed $N$ and $\ell$.

\begin{lemma}\label{lemmares} (a) Each $g\in O(N)$ for which $\chi_\mu(g)=0$ for all boundary
diagrams $\mu$ of $\LNell$ defines a trace on $\qBb_\infty$ determined
by $tr(p_\mu)=\chi_\mu(g)/\chi_{[1]}(g)^n$, where $p_\mu$ is a minimal
projection of $\qBb_{n,\mu}$.

(b) For given $\la\in\tLNell_n$ the coefficients $\bNell$ are uniquely determined by the equations
$\chi^{U(N)}_\la(g)=\sum_\mu\bNell\chi_\mu(g)$
for all $g$ as in (a), where the summation goes over
all diagrams $\mu$ in $\LNell$ with $n, n-2,\ ... $ boxes.
\end{lemma}
\
$Proof.$ The formula in statement (a) determines a trace on $\qBb_n$
for each $n$.
To show that these formulas are compatible with the standard embeddings
we observe that a minimal idempotent $p_\mu\in\qBb_{n,\mu}$
is the sum of minimal idempotents $e_\la\in\qBb_{n+1,\la}$
where $\la$ runs through all diagrams in $\LNell$ obtained by adding or
removing a box to/from $\la$, see Eq \ref{decomposition}
and the remarks below that theorem.
Evaluating the traces of these idempotents
and multiplying everything by $\chi_\la(g)^{n+1}$, equality of the
traces is equivalent to
$$\chi_\mu(g)\chi_{[1]}(g)=\sum_\la \chi_\la(g).$$
By the usual tensor product rule for orthogonal groups,
the left hand side would be equal to the sum of characters
corresponding to $all$ diagrams $\la$ which differ from $\mu$ by only
one box. It is easy to check that this differs from the sum above only
by boundary diagrams, for which the characters at $g$ is equal to 0.
This shows (a).

For (b), we first show that $tr(p_\la)=\chi_\la^{U(N)}(g)/\chi^{U(N)}_{[1]}(g)^n$
for $p_\la\in\bar H_{n,\la}$ a minimal idempotent and $tr$ a trace as in (a).
 As the weight vector for $\qBb_{n+2N}$ is
a multiple of the one of $\qBb_n$, for $n$ large enough, the same
must also hold for the weight vectors of $\bar H_n$ and $\bar H_{n+2N}$,
by periodicity of the inclusions. Hence these weight vectors must
be eigenvectors of the inclusion matrix for 
$\bar H_n\subset\bar H_{n+2N}$. As this  inclusion matrix is
just a block of the $2N$-th  power of the fusion matrix
of the vector representation for the corresponding type $A$
fusion category, its entries must be given by $U(N)$ characters
of a suitable group element. 
To identify these elements, it suffices to observe that
 the antisymmetrizations of the vector representation, labeled by the Young diagrams
 $\la =[1^j]$, $1\leq j\leq N$, remain irreducible as $O(N)$ modules.
This means the corresponding Hecke algebra idempotent remains
a minimal idempotent also in $\qBb_j$. Hence $tr(p_{\la})=\chi_{\la}^{U(N)}(g)$
for $\la=[1^j]$ and $1\leq j\leq N$.
 But as the antisymmetrizations
generate the representation ring of $U(N)$, and also of the corresponding
fusion ring, the claim follows for general $\la$. For more details, see e.g.
\cite{GW}

Recall that the coefficient $\bNell$ can be defined as the rank of $p_\la$
in an irreducible $\qBb_{n,\mu}$ representation. So obviously the formula
in the statement holds for any $g$ as in (a). Examples for such $g$ come
from $exp(\x)$ with $\x\in M^*=\ell^{-1}Q$ for which the character
is given by the expression $\chi_\la(\x)$ as in  Section \ref{traceweights}.
As the columns of the orthogonal $S$-matrix are linearly independent,
this would identify $SO(N)$ representations. If $N$ is odd, the two $O(N)$
representations which reduce to the same $SO(N)$ representation
are labeled by Young diagrams with opposite parities. Hence only one
of them can occur in the decomposition of a given $U(N)$ representation.
A similar argument also works in the symplectic case.

For $N$ even, we can have two diagrams $\la$ and $\la^\dagger$
with the same $SO(N)$ character, where one of them, say $\la$ has less
than $k$ rows. They can be distinguished by elements $g\in O(N)\backslash{}
SO(N)$ for which $\chi_{\la^\dagger}(g)=-\chi_\la(g)$.
 It is well-known that such elements $g$ must have eigenvalues $\pm 1$,
and $\chi_\la(g)$ is given
by the character formula for $Sp(2k-2)$ in the remaining $2k-2$ eigenvalues
(see  \cite{Wy}). It follows from the invertibility of the $S$-matrix
for $Sp(2k-2)$ at level $\ell/2-k$ (see \cite{Kc}) that we can
identify those diagrams $\la$ by evaluating $\chi_\la^{Sp(2k-2)}(\x/\ell)$ for
 $\x\in\Z^{k-1}$ with $\ell/2>x_1>x_2>\ ...\ >x_{k-1}>0$, and that those
elements satisfy the boundary condition $\chi^{Sp(2k-2)}_\la(\x)=0$
for any boundary diagram $\la$.

\medskip

The lemma above is illustrated in the following section for a number of
explicit examples. We can also give a closed formula for the restriction 
coefficients, using a well-known quotient map for fusion rings (even though
in our case, the quotient ring does not correspond to a tensor category
as far as we know). In  the context of fusion rings, this is known as
the Kac-Walton formula; for type $A$ see also e.g. \cite{GW}.
In our case, we need to use a slightly different affine reflection group $\W$.
In the orthogonal case $N=2k$ and $N=2k+1$ it is given by the semidirect product
of $\ell \Z^k$ with the Weyl group of type $B_k$.
In the symplectic case, it is given by the semidirect
product of $\ell Q$ with the  Weyl group of type $B_k$. 
As usual, we define the dot action of $\W$ on $\R^k$
by $w.\x=w(\x+\rho)-\rho$, where $\rho$ is half the sum
of the positive roots of the corresponding Lie algebra,
with the roots embedded into $\R^k$ as described above,
and $\ve$ is the usual sign function for reflection groups.
This can be extended to an action on the labeling set of
$O(N)$ representations by identifying a Young diagrams with
$\leq k$ rows with the corresponding vector in $\Z^k$,
and by using the restriction rule from $O(N)$ to $SO(N)$
in the other cases. See also e.g. Lemma 1.7 in \cite{Wquot}
for more details.

\begin{theorem}\label{result} With notations as above,
the restriction multiplicity $\bNell$ for $N=2k+1>0$ and
$N=-2k$  is given by
$$\bNell=\sum_{w\in \W} \ve(w)b^\la_{w.\mu}(N).$$
If $N=2k>0$, we have to replace $\ve$ by $\vet$
(see Section \ref{anothers}) in
the formula above.
\end{theorem}

$Proof.$  Looking at the character formulas, we see that
an action of an element $w$ of the finite reflection group
 on $\la$ just changes the character
by the sign of $w$. Moreover, by definition of the elements
$\x$ we have that $\chi_\la(\x)=\chi_{\la+\mu}(\x)$ for any $\mu\in M$.
It follows that  $\chi_{w.\la}(\x)=\ve(w) \chi_\la(\x)$ for
all $\x\in M^*$ and $w\in\W$. Hence summing over the $\W$-orbits, we obtain
for any $\x\in M^*$, $\la\in\tLNell$ and $\mu\in\LNell$ that
$$\chi^{U(N)}_\la(\x)=\sum_{\gamma} b^\la_\gamma(N)\chi_\gamma
= \sum_\mu (\sum_w b^\la_{w.\mu}(N))\chi_\mu.$$
The claim now follows from this and Lemma \ref{lemmares}.

\section{Examples and other approaches}

\subsection{ The case $N=2$:} 
This corresponds to the Goodman-de la Harpe-Jones subfactors
for type $D_{\ell/2+1}$, 
where $\ell>2$ has to be even. It follows from
our theorem that the even vertices of the principal
graph are labeled by the Young diagrams $\la$ with an even number $n$
of boxes, at most two rows and with $\la_1-\la_2\leq \ell-2$;
there are $(\ell-2)/2$ such diagrams.
Their dimensions are given by $\tilde d_k=[2k+1]$, $0\leq k<(\ell-2)/2$.

Moreover, one checks that $\Lambda(2,\ell)$ consists of Young
diagrams $[j]$ with $0\leq j\leq \ell/2-1$ and of $[1^2]$, one column with 2 boxes,
with dimensions $d_{[j]}=2\cos j\pi/\ell$ for $j>0$ and dimension equal to 1 for the
remaining cases (i.e. for $\emptyset$ and for $[1^2]$).
The restriction rule (i.e. principal graph) follows from writing
the dimensions as
$$\tilde d_k = 2 \cos \tilde k \pi/\ell \ +\ 2\cos (\tilde k-2)\pi/\ell +\ ...\ +1,$$
where $\tilde k= $ min $\{ k, \ell/2-k\}$.  Indeed, this determines the graph
completely except for whether to pick the diagram $\emptyset$ or $[1^2]$
for the object with dimension 1. It follows from the restriction rule
$O(2)\subset U(2)$ that we take $\emptyset$  for $j$ even, and  $[1^2]$
for $j$ odd. To calculate the index one can check  by elementary means that
$\sum_{\la \ even} d_\la^2 = \ell/2$. Moreover, it is well-known that
the  sum $\sum_{\la \ even} \tilde d_\la^2$  over even partitions for $\sl_2$ is equal to
$\ell/4\sin^2\pi/\ell$.  Hence we obtain as index
$[\M:\N]=1/2\sin^2\pi/\ell$.

\subsection{ The case $N=3$} 
\begin{figure}
\centering
\includegraphics[scale=.7]{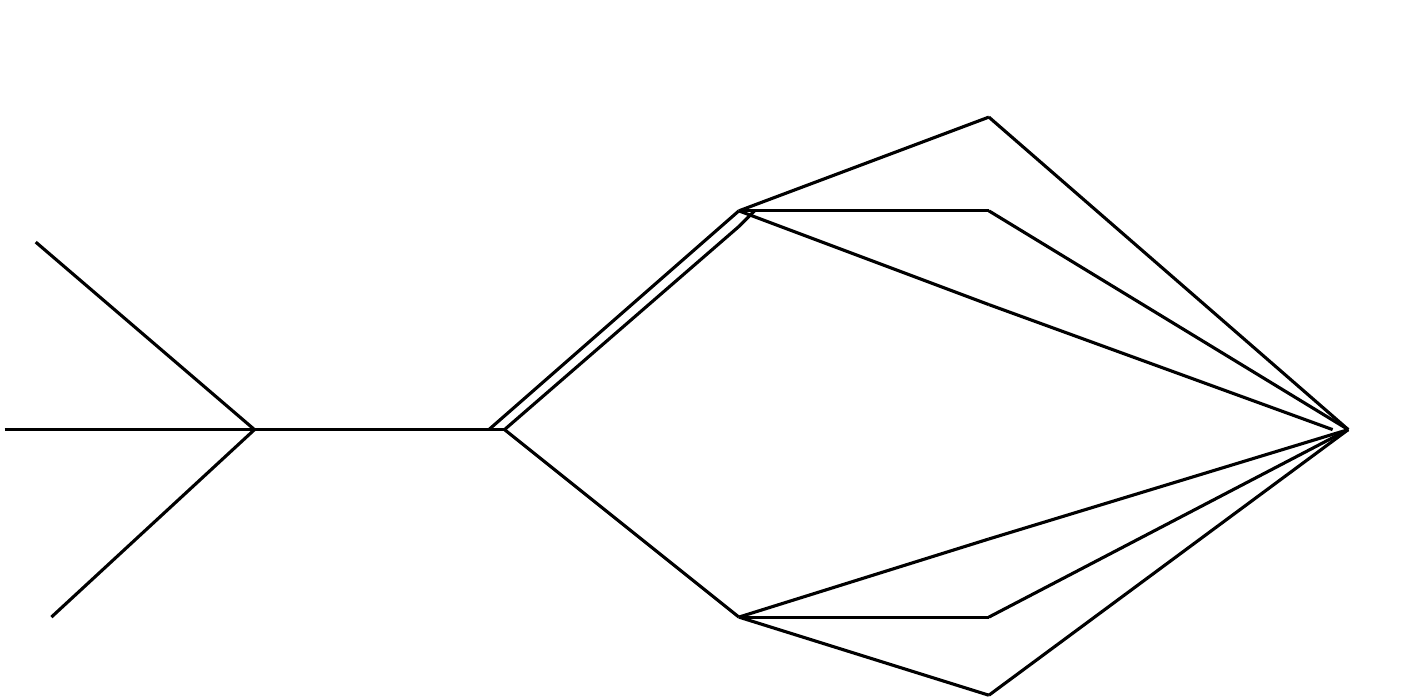}
\caption{SO(3) for $\ell = 9$}
\label{fig:SO(3)}
\end{figure}

It is also fairly elementary to work out this case in detail.
Recall that by Weyl's dimension formula we have
$$\tilde d_\la = \frac{ [\la_1-\la_2+1][\la_2-\la_3+1][\la_1-\la_3+2]}{[1]^2[2]}.$$
Now observe that the product of two $q$-numbers is given by the tensor product
rules for $\sl_2$, i.e. we have for $n\geq m$ that 
$[n][m]=[n+m-1] + [n+m-3] + ... + [n-m+1]$.
As an example, we have
$$\tilde d_{[4]}\ =\ \frac{[6][5]}{[2]}\ =\ \frac{[10]+[8]+[6]+[4]+[2]}{[2]}\ =\ [9]+[5]+[1],$$
i.e. the fourth antisymmetrization of the vector representation of $U(3)$ decomposes
as a direct sum of the one-, five- and nine-dimensional representation of $SO(3)$.
One similarly can show the well-known result that the adjoint representation of $SU(3)$,
labeled by the Young diagram $[2,1]$ decomposes into the direct sum of the three- and
the five-dimensional representation of $SO(3)$, i.e. $\p$ is the five dimensional representation
of $SO(3)$. Hence we get from Theorem \ref{explicitindex} that the index is equal to
$$[\M:\N]\ =\ \frac{\ell}{4^2\sin^2 (2\pi/\ell)\ \sin^2 (\pi/\ell)}.$$
We note that here as well as in the other examples, the dimensions (i.e. entries of
the Perron-Frobenius vectors) are given by $|\tilde d_\la|$ for even vertices,
and by $\sqrt{[\M:\N]}|d_\mu|$ for odd vertices, with $\tilde d_\la$ and $d_\mu$ as in
Lemma \ref{weightcharacters}.
To consider explicit examples, the first nontrivial case for $N=3$ occurs for $\ell=7$.
We leave it to the reader to check that in this case the first principal graph is given
by the Dynkin graph $D_8$. A more interesting graph is obtained for $\ell=9$, see
Fig \ref{fig:SO(3)}. Here we have the three invertible objects of the $SU(3)_6$
fusion category, including the trivial object (often denoted as $*$) on the left;
they generate a group isomorphic to $\Z/3$.
The vertices with the double edge are labeled by the object corresponding to
the 5-dimensional representation of $SO(3)$ and the diagram $[4,2]$ for $SU(3)_6$.
This is the only fixed point under the $\Z/3$ action given by the invertible objects
(or, in physics language, the currents).  It would be interesting to see whether
one can carry out an orbifold construction in this context related to the one
in \cite{EK1}.

\subsection{The case $N= 4$} 
\begin{figure}
\centering
\includegraphics[scale=.7]{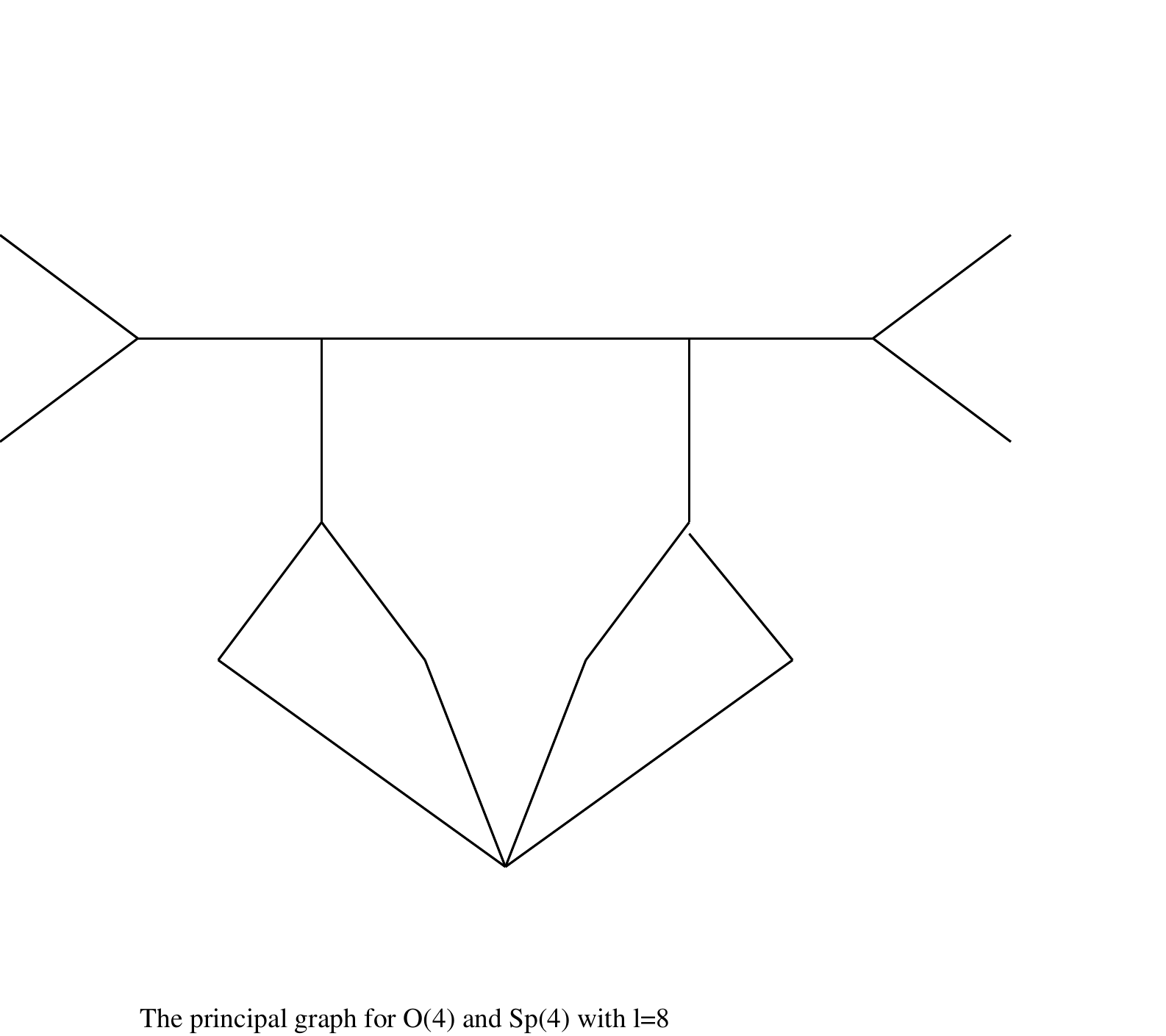}
\caption{O(4) and Sp(4) for $\ell = 8$}
\label{fig:O(4)}
\end{figure}

This is the first nontrivial and, apparently, new case
corresponding to even-dimensional orthogonal groups. As we shall see, somewhat 
surprisingly, the corresponding construction for $SO(4)$ does not seem to work.
We do the case with $\ell=8$ in explicit detail. It is not hard to check that we already
get the periodic inclusion matrix for $n=12$. As we consider an analog of the
restriction to $O(4)$ for which the determinant can be $\pm 1$, we should,
strictly speaking, consider a fusion category for $SU(4)\times \{\pm 1\}$.
We shall actually use the Young diagram notation for representations of
$U(4)$. For $n=12$ we have the invertible objects labeled by
$[3^4]$, $[4^3]$, $[5^21^2]$ and $[62^3]$ (i.e. e.g. the last diagram has six boxes
in the first and two  boxes in the second, third and fourth row). They generate
a subgroup isomorphic to $\Z/4$.
It follows from the 
$O(4)$ restriction rules that $[3^4]$ and $[5^21^2]$ contain the determinant representation,
 and $[4^3]$ and $[62^3]$ contain the trivial representation
 as one-dimensional $O(4)$  subrepresentations. This allows us to calculate the
restrictions for representations of each $\Z/4$ orbit simultaneously. As usually for
at least one element of each orbit the ordinary restriction rules still hold, it makes
the general calculations easier. The principal graph can be seen in Fig. \ref{fig:O(4)}.
As in the $N=3$ example, the one-dimensional currents, including the trivial object
$*$ appear as the left- and right-most vertices in the graph. The lowest vertex
corresponds to the $O(4)$-object $[2]$ which is connected to the objects
in the $\Z/4$-orbit $\{[2,1^2], [3,1], [4,3,1], [3,3,2]\}$. We also note that
we get the same graph for the $Sp(4)$ case $N=-4$ for $\ell=8$. However,
for other roots of unities, already the indices of the subfactors differ which
are given by
$$O(4): \ \frac{2\ell}{4\sin^2(3\pi/\ell)\ 4\sin^2 (2\pi/\ell)\ 16\sin^4 (\pi/\ell)}\hskip 3em
Sp(4):\ \frac{\ell}{4\sin^2 (2\pi/\ell)\ 4\sin^2 (\pi/\ell)}.$$

It was originally thought that we should also be able to get fusion category analogs for the
restriction from $SU(N)$ to $SO(N)$ for $N$ even. It is easy to check that this is not
possible for $O(2)$. Some initial checks also seem to suggest a similar phenomenon for
higher ranks. E.g. using the same element $\rhoc$ in the $SO(N)$ character formula
would give dimension functions which are not invariant under the $D_N$ diagram
automorphism.

\subsection{Related results} We discuss several results related to our findings.
Our original motivation was to construct subfactors related to twisted loop groups.
It was shown in R. Verill's  PhD thesis \cite{V} that it is not possible to
construct a fusion tensor product for representations of twisted loop groups.
However, it seemed reasonable to expect that representations of twisted
loop groups could become a module category over representations of
their untwisted counterparts. Many results, in particular about the
combinatorics of such categories can be found in the context
of boundary conformal field theory in papers by Gaberdiel, Gannon,
Fuchs, Schweigert, diFrancesco, Petkova, Zuber and others
(see e.g. \cite{GG}, \cite{FS}, \cite{PZ} and the papers cited therein).

In the mathematics literature, one can find closely related results 
in the papers \cite{X} and \cite{Wa}. Here the authors construct
module categories via a completely different
approach  in the context of loop groups. E.g. the formulas
at the end of \cite{X} for the special case $N=3$ differ only by a factor 3 
(which can be explained, see Remark \ref{indexremark}) from our formulas for $N=3$ for even level
(together with Corollary \ref{localind}),
 modulo misprints; similar formulas for the symplectic case as well
as restriction coefficients also
appear at the end of \cite{Wa}.
We can not get results corresponding to the odd level cases in \cite{X}.
The combinatorics there suggests that this would require considering
an embedding of $Sp(N-1)$ into $SU(N)$ under which the vector
representation would not remain irreducible. In contrast,
we can also construct module categories for $\ell-N$ odd, which
would correspond to odd level; however, these categories are
not  unitarizable (which follows from Lemma \ref{weightpositivity})
and they have different fusion rules.
However, we do get fairly general formulas for the index and principal
graphs of this type of subfactors in the unitary case, 
which was one of the problems posed
in \cite{X}. These formulas were known to this author as well as to
Antony Wassermann at least back in 2008 when they had discussions
about their respective works in Oberwolfach and at the Schr\" odinger
Institute.

We  close this section by mentioning that while our results for
$N>0$ odd and $N<0$ even are in many ways parallel to results
obtained via other approaches in connection with twisted loop 
groups, there does not seem to be an obvious analog for our
results for $N>0$ even. E.g. the combinatorial results
in \cite{GG} for that case seem to be different to ours.

\subsection{Conclusions and further explorations} We have constructed module
categories of fusion categories of type $A$ via deformations of centralizer
algebras of certain subgroups of unitary groups. We have also classified
when they are unitarizable, and we have constructed the corresponding
subfactors.  These deformations
are compatible with the Drinfeld-Jimbo deformation of the unitary group
but not with the Drinfeld-Jimbo deformation of the subgroup.
Most of the deformation was already 
done in \cite{Wqalg} via elementary methods. In principle, at least,
it should be possible to use this elementary approach also for
other inclusions. However, this might become increasingly tedious.

As we have seen already  in Section \ref{Molev}, it should be possible
to get a somewhat more conceptual approach using different deformations
of the subgroup, see \cite{N}, \cite{Mo},  \cite{L1}, \cite{L}, \cite{IK} and 
references therein. In particular in the work of Letzter,
such deformations  via co-ideal algebras have been defined for 
a large class of embeddings of a semisimple Lie algebra into another one.
At this point, it does not seem obvious how to define $C^*$-structures
in this setting, and additional complications arise as these coideal algebras
are not expected to be semisimple at roots of unity.
Nevertheless, the results in this and other  papers such as
\cite{X}, \cite{Wa} would seem to suggest
that similar constructions might be possible also in a more general setting.
It would also be interesting to see whether this could lead
to a general solution of Problem 5.5 in \cite{O}.

\bibliographystyle{plain}

\begin{thebibliography}{99}
\bibitem[BW] {BW} Birman, J. and Wenzl, H., Braids, link polynomials and
a new algebra, Trans. Amer. Math. Soc. 313(1) (1989) 249-273
\bibitem[Bi]{Bi}  Bisch, D., Bimodules, higher relative commutants and the fusion algebra associated
to a subfactor,  pp. 13–63 in Operator algebras and their applications (Waterloo, ON, 1994/1995),
edited by P. A. Fillmore and J. A. Mingo, Fields Inst. Commun. 13, Amer. Math. Soc., Providence,
RI, 1997.
\bibitem[Bl] {Bl}  Blanchet, Ch., Hecke algebras, modular
categories and $3$-manifolds quantum invariants.
Topology  39  (2000),  no. 1, 193--223.
\bibitem[BEK]{BEK} B\" ockenhauer, Jens; Evans, David E.; Kawahigashi, Yasuyuki 
On $\alpha$-induction, chiral generators and modular invariants for subfactors. 
Comm. Math. Phys. 208 (1999), no. 2, 429–487.
\bibitem[Br] {Brauer} Brauer, R., On algebras which are connected with the
semisimple continuous groups, Ann. of Math. {\bf 63} (1937), 854-872.
\bibitem[E]{E} Erlijman, Juliana New subfactors from braid group representations. 
Trans. Amer. Math. Soc. 350 (1998), no. 1, 185–211.
\bibitem[EW]{EW} Erlijman, Juliana; Wenzl, Hans, Subfactors from braided 
$C^*$C∗ tensor categories. Pacific J. Math. 231 (2007), no. 2, 361–399
\bibitem[EK1]{EK1} Evans, David E.; Kawahigashi, Yasuyuki Orbifold subfactors from Hecke algebras. Comm. Math. Phys. 165 (1994), no. 3, 445–484.
\bibitem[EK]{EK} Evans, D.E., Kawahigashi, Y.: Quantum symmetries on operator algebras.
 Oxford: Oxford University Press, 1998
\bibitem[FS]{FS}  Fuchs, J.,  Schweigert, C., Solitonic sectors, alpha-induction and symmetry breaking
boundaries, Phys. Lett. B490 (2000) 163
\bibitem[GG]{GG} Gaberdiel, M., Gannon, T., Boundary states for WZW models. 
Nucl.Phys. B639 (2002) 471-501
\bibitem[GHJ]{GHJ} Goodman, Frederick M.; de la Harpe, Pierre; Jones, Vaughan F. R. 
Coxeter graphs and towers of algebras. Mathematical Sciences Research Institute 
Publications, 14. Springer-Verlag, New York, 1989. 
\bibitem[GW]{GW} Goodman, Frederick M.; Wenzl, Hans Littlewood-Richardson coefficients for Hecke algebras at roots of unity. Adv. Math. 82 (1990), no. 2, 244–265.
\bibitem[IK]{IK} Iorgov, N. Z.; Klimyk, A. U. Classification theorem on irreducible 
representations of the $q$-deformed algebra $U_q'({\rm so}_n)$ Int. J. Math. Math. Sci. 2005, no. 2, 225–262.
\bibitem[J]{Jo} Jones, V.F.R, Index for subfactors, Invent. Math. 72, (1983) 1-25
\bibitem[Kc]{Kc} V. Kac, Infinite dimensional Lie algebras, 3rd edition,
Cambridge University Press.
\bibitem[Ks]{Kassel} Kassel, Ch., Quantum groups, Springer
\bibitem[Ko]{Ko} Koike, K., Principal specializations of the classical
groups and $q$-analogs of the dimension formulas.  Adv. Math.  125  (1997),
no. 2, 236--274.
\bibitem[KT]{KT} Koike, K. and Terada, I., Young-diagrammatic methods for
the representation theory of the classical groups of type
$B\sb n,\;C\sb n,\;D\sb n$.  J. Algebra  107  (1987),  no. 2, 466--511.
\bibitem[L1]{L1} Letzter, G., Subalgebras which appear in quantum Iwasawa
decompositions, Can. J. Math. Vol. 49 (6), 1997 pp. 1206-1223.
\bibitem[L2]{L} Letzter, G., Coideal subalgebras and quantum symmetric pairs.
New directions in Hopf algebras,  117--165, Math. Sci. Res. Inst. Publ., 43.
\bibitem[LR]{LR} Longo, R., Roberts, J., A theory of dimension, $K$-theory, 11 (1997) 103-159
\bibitem[Mo]{Mo} Molev, A.I., A new quantum analog of the Brauer algebra. 
 Czech. Journal of Physics 53 (2003), 1073-1078
\bibitem[N]{N} Noumi, N.,  Macdonald’s symmetric polynomials as zonal spherical functions on
quantum homogeneous spaces, Adv. Math. 123 (1996), 16–77.
\bibitem[O]{O} Ostrik, Victor Module categories, weak Hopf algebras and modular invariants.
 Transform. Groups 8 (2003), no. 2, 177–206
\bibitem[PZ]{PZ} Petkova, V., Zuber, J.-B., Conformal field theories, graphs and quantum algebras,
hep-th/0108236
\bibitem[RW]{RW} Ram, A. and Wenzl, H., \emph{Matrix units for centralizer
algebras} J. Algebra. 102 (1992), 378-395.
\bibitem[Su]{Su} Sunder, V. S. A model for AF algebras and a representation of the Jones projections. J. Operator Theory 18 (1987), no. 2, 289–301. 
\bibitem[TuW]{TuW} Tuba, I and Wenzl, H., On braided tensor categories of
type $BCD$, J. reine angew. Math. 581 (2005), 31-69.
\bibitem[Tu]{Turaev} Turaev, V. Quantum invariants, DeGruyter
\bibitem[V]{V} Verril, R., Positive energy representations of $L^\sigma SU(2r)$ amd orbifold fusion,
Pd.D. Thesis, University of Cambridge, 2001
\bibitem[Wa]{Wa} Wassermann, A.,  Subfactors and Connes fusion for twisted loop groups, arXiv:1003.2292
\bibitem[W1]{WHe} Wenzl, H., \emph{Hecke algebras of type $A_n$ and subfactors},
Invent. Math 92, 349-383 (1988).
\bibitem[W2]{w1} Wenzl, H., \emph{Quantum Groups and Subfactors of type $B$,
$C$, and $D$}, Comm. Math. Phys. 133, 383-432 (1990).
\bibitem[W3]{w2} Wenzl, H., \emph{On the structure of Brauer's centralizer
algebras}, Ann. of Math., 128, 173-193 (1988).
\bibitem[W4]{Wqalg} Wenzl, H., A $q$-Brauer algebra, arXiv:1102.3892
\bibitem[W5]{Wquot} Wenzl, H., Quotients of representation rings,
Represent. Theory 15 (2011), 385-406. 
arXiv:1101.5887
\bibitem[Wy]{Wy} Weyl, H., The classical groups, Princeton University Press.
\bibitem[X1]{X1}  Xu, F.,  New braided endomorphisms from conformal inclusions. Commun. Math. Phys. 192,
347–403 (1998)
\bibitem[X]{X} Xu, Feng On affine orbifold nets associated with outer automorphisms. 
Comm. Math. Phys. 291 (2009), no. 3, 845–861.
\end{thebibliography}

\end{document}